# A SHRINKING LEMMA FOR INDEXED LANGUAGES

Robert H. Gilman

ABSTRACT. This article presents a lemma in the spirit of the pumping lemma for indexed languages but easier to employ.

## Section 1. Introduction.

The pumping lemma for context-free languages has been extended to stack languages [O] and indexed languages [H], but these generalizations are rather complicated. In this article we take a slightly different approach by concentrating only on that part of the context–free pumping lemma which says that if $uvwxy \in L$, then $uwy \in L$, and by employing a theorem on divisibility of words which is not used in [O] or [H]. Our result, Theorem A, is relatively easy to state and strong enough to verify the examples given in [H] of languages which are not indexed. On the other hand it does not afford a proof that the finiteness problem for indexed languages is solvable as does [H, Theorem 5.1].

Indexed languages were introduced by Aho [A1], [A2]. A brief introduction appears in [HU, Chapter 14]. Our original motivation for Theorem 1 was the investigation of finitely generated groups for which the language of words defining the identity is indexed.

## Section 2. A Result on Indexed Languages.

Before stating our result we fix some notation. $\Sigma$ is a finite aphabet, $|w|$ is the length of $w \in \Sigma^*$, and for each $a \in \Sigma$, $|w|_a$ is the number of $a$'s in $w$.

**Theorem A.** *Let $L$ be an indexed language over $\Sigma$ and $m$ a positive integer. There is a constant $k > 0$ such that each word $w \in L$ with $|w| \geq k$ can be written as a product $w = w_1 \cdots w_r$ for which the following conditions hold.*

(1) $m < r \leq k$.
(2) *The factors $w_i$ are nonempty words.*

---

The author was partially supported by NSF grant DMS-9401090.

Typeset by $\mathcal{AMS}$-TeX

  (3) *Each choice of $m$ factors is included in a proper subproduct which lies in $L$.*

By (3) we mean that the chosen factors occur in a product $w_{i_1} \cdots w_{i_t} \in L$ with $1 \leq i_1 < \ldots < i_t \leq r$ and $m \leq t < r$. The proof of Theorem A is given in the next section.

**Corollary 1.** *Let $L$ be an indexed language. There is a constant $k > 0$ such that if $w \in L$ and $|w| > k$, then there exists $v \in L$ with $(1/k)|w| \leq |v| < |w|$.*

*Proof.* Take $m = 1$ in Theorem A and choose a factor of maximum length. □

By taking $m$ to be the number of letters in $\Sigma$ and arguing similarly we obtain a result on the Parikh mapping.

**Corollary 2.** *Let $L$ be an indexed language over $\Sigma$. There is a constant $k > 0$ such that if $w \in L$ and $|w| > k$, then there exists $v \in L$ with $(1/k)|w|_a \leq |v|_a \leq |w|_a$ for each $a \in \Sigma$ and $|v|_a < |w|_a$ for some $a \in \Sigma$.*

Corollary 1 has the following immediate consequence.

**Corollary 3.** [H, Theorem 5.2] *If $f$ is a strictly increasing function on the positive integers, and $L = \{a^{f(n)}\}$ is an indexed language, then $f = O(k^n)$ for some positive integer $k$.*

**Corollary 4.** [H, Theorem 5.3] *The language $L = \{(ab^n)^n \mid n \geq 1\}$ is not indexed.*

*Proof.* Suppose $L$ is indexed, and apply Theorem A to $L$ with $m = 1$. Pick $w = (a^n b)^n$ with $n > k$ and consider the decomposition $w = w_1 \cdots w_r$. As $r \leq k$, at least one factor $w_i$ must contain two or more $a$'s. Choose that $w_i$ to be in the proper subproduct $v$. But then $v$ contains a subword $ab^n a$, which is impossible as $v \neq w$. □

## Section 3. Proof.

The proof of Theorem A depends on a result about divisibility of words. We say that $v$ divides $w$ and write $v \prec w$ if $v$ is a subsequence of $w$. For example $ac \prec abc$. By a theorem of Higman [SS, Theorem 6.1.2] every set of words defined over a finite alphabet and pairwise incomparable with respect to divisibility is finite. We will use this result in the following form.

**Lemma 1.** *Let $m$ be a positive integer and $Y$ a language over a finite alphabet $\Delta$. $Y$ contains a finite subset $X$ with the property that for any* $y \in Y - X$ with $m$ letters distinguished there is an $x \in X$ such that $x \npreceq y$ and $x$ includes all the distinguished letters of $y$.

*Proof.* Let $\Delta'$ be the union of $\Delta$ with $m$ pairwise disjoint copies of itself, and define $Y'$ be the language of all words over $\Delta'$ which project to $Y$ and contain exactly one letter from each of the $m$ copies of $\Delta$. By Higman's theorem $X'$, the set of all words in $Y'$ each of which is not divisible by any word in $Y'$ except itself, is finite. For any $y' \in Y'$ if we take $x'$ to be a word of minimum length among all words in $Y'$ dividing $y'$, then $x' \in X'$. Further $x'$ contains all the letters of $y'$ from $\Delta' - \Delta$.

Define $X$ to be the union of the projection of $X'$ to $\Delta^*$ with the set of all words in $Y$ of length less than $m$. Suppose that $y \in Y - X$ has $m$ distinguished letters. Since $|y| \geq m$, we can pick $y' \in Y'$ projecting to $y$ so that the distinguished letters of $y$ correspond to the letters of $y'$ in $\Delta' - \Delta$. By the preceding paragraph $y'$ is divisible by an $x' \in X'$ which contains those letters. It follows that the projection of $x'$ to $\Sigma^*$ is the desired word $x$. □

Notice that $x$ might be a subsequence of $y$ in more than one way. Lemma 1 asserts only that there is some subsequence of $y$ which includes the distinguished letters and whose product is $x$.

Fix an indexed language $L$ over $\Sigma$, and let $G$ be an indexed grammar for $L$. Let $G$ have sentence symbol $S$, nonterminals $N$, and indices $F$. $(NF^* + \Sigma)^*$ is the set of sentential forms. By [A1, Theorem 4.5] we may assume $G$ is in normal form, i.e.,

(1) $S$ does not appear on the righthand side on any production;
(2) There are no $\epsilon$-productions except perhaps $S \to \epsilon$;
(3) Each production has one of the forms $A \to BC$, $Af \to B$, $A \to Bf$, or $A \to a$, where $A, B, C \in N$, $f \in F$, and $a \in \Sigma$.

We are using the definition of indexed grammar from [HU]; this definition is slightly different from the original.

We write $\alpha \xrightarrow{*} \beta$ to indicate that the sentential form $\beta$ can be derived from the sentential form $\alpha$ via productions of $G$, and we use $\beta \cdot \omega$ to denote the sentential form obtained by appending the index string $\omega$ to the index string of every nonterminal in the sentential form $\beta$. It follows from the way derivations are defined in indexed grammars that if $\alpha \xrightarrow{*} \beta$, then $\alpha \cdot \omega \xrightarrow{*} \beta \cdot \omega$. Conversely if $\alpha \cdot \omega \xrightarrow{*} \beta \cdot \omega$ and if every nonterminal occurring in that derivation has an index string with suffix $\omega$, then $\alpha \xrightarrow{*} \beta$.

**Lemma 2.** *Let $m$ be a positive integer and $A\omega$ a sentential form in $NF^*$. There is a finite set of sentential forms $X \subset (N + \Sigma)^*$ with the property that* *if $A\omega \xrightarrow{*} \beta \in (N+\Sigma)^* - X$, and $m$ symbols of $\beta$ are distinguished, then there is $\alpha \in X$ such that $A\omega \xrightarrow{*} \alpha \precneqq \beta$, and $\alpha$ includes all the distinguished symbols of $\beta$.*

*Proof.* Apply Lemma 1 to the language of all sentential forms in $(N+\Sigma)^*$ derivable from $A\omega$. □

Consider a derivation $S \xrightarrow{*} w \in L$, and let $\Gamma$ be the corresponding derivation tree. Let each vertex $p$ of $\Gamma$ have label $\lambda(p)$, and define a subtree $\Gamma(p)$ with root $p$ as follows. If $\lambda(p)$ is a terminal or nonterminal, then $\Gamma(p)$ consists of $p$ and all its descendants. Otherwise $\lambda(p) = Af\omega$ for some nonterminal $A$, index $f$, and string of indices $\omega$. In this case along each path emanating from $p$ there will be a first vertex, perhaps a leaf of $\Gamma$, at which $f$ is consumed. Define $\Gamma(p)$ to be the union of all the paths from $p$ up to and including these first vertices. The subtrees $\Gamma(p)$ play an important role in [H]; we will use them here in a slightly different way than they are used there.

Let $\gamma(p)$ be the sentential form obtained by concatenating the labels of the leaves of $\Gamma(p)$ in order; if $p$ is a leaf, $\gamma(p) = \lambda(p)$. Since $\Gamma(p)$ is a subtree of a derivation tree, $\lambda(p) \xrightarrow{*} \gamma(p)$. If $\lambda(p) = Af\omega$, then by construction all vertices of $\Gamma(p)$ except its leaves have labels of the form $B\omega' f\omega$. The leaves are labelled by terminals or labels form $B\omega$. Deleting all the suffixes $\omega$ yields a derivation tree for $Af \xrightarrow{*} \beta(p)$ where $\gamma(p) = \beta(p) \cdot \omega$. Extend the definition of $\beta(p)$ to all other vertices $p$ of $\Gamma$ by defining $\beta(p) = \gamma(p)$ when $\lambda(p)$ is a terminal or nonterminal.

It follows from Lemma 2 that there is a finite set of sentential forms $Z \subset (N+\Sigma)^*$ such that for any of the finitely many sentential form $A\omega \in N \cup NF$ if $A\omega \xrightarrow{*} \beta \in (N+\Sigma)^* - Z$ and $m$ symbols of $\beta$ are distinguished, then there is $\alpha \in Z$ such that $A\omega \xrightarrow{*} \alpha \precneqq \beta$, and $\alpha$ includes all the distinguished symbols of $\beta$. Since it does no harm to enlarge $Z$, we may assume $Z$ contains all elements of $(N+\Sigma)^*$ of length at most $m$.

**Lemma 3.** *Let $C \geq 2$ be an upper bound for the lengths of elements of $Z$. Suppose $\beta(p) \notin Z$ but $\beta(q) \in Z$ for all vertices $q$ which are proper descendants of $p$, then $|\beta(p)| \leq C^2$.*

*Proof.* If $p$ is a leaf, then $|\beta(p)| = 1$. Suppose $p$ has two descendants, $q_1, q_2$. It follows from the normal form for $G$ that $\beta(p) = \beta(q_1)\beta(q_2)$, and consequently $|\beta(p)| \leq 2C$. Finally if $p$ has a single descendant, $q$, then the derivation $\lambda(p) \xrightarrow{*} \gamma(p)$ begins with application of a production of the form $A \to a$, $Af \to B$ or $A \to Bf$. In the first case $|\beta(p)| = |a| = 1$. In the second case $\lambda(p)$ must be $Af\omega$ whence $\beta(p) = B$ and again $|\beta(p)| = 1$.

Consider the last case. We have $\lambda(p) = A\omega$ and $\lambda(q) = Bf\omega$. Further $\beta(p)$ is the product of the terms $\beta(q')$ as $q'$ ranges over the leaves of $\Gamma(q)$. Since $\beta(q) \in Z$, there are at most $C$ terms; and as each $\beta(q') \in Z$, we have $|\beta(p)| \leq C^2$. □

To complete the proof of Theorem A choose $k = C^2 + 2$ and suppose $S \overset{*}{\to} w \in L$ with $|w| \geq k$. Let $\Gamma$ be the corresponding derivation tree and $p_0$ its root. Clearly $\beta(p_0) = w \notin Z$, and so we may choose $p$ to satisfy the hypothesis of Lemma 3. Note that $\beta(p) \notin Z$ implies $|\beta(p)| > m$; in particular $p$ is not a leaf.

If $\lambda(p) = A$, then $\beta(p) = a_1 \cdots a_t$ is a subword of $w$ and $m < t \leq C^2$. Consequently $w = w'a_1 \cdots a_t w''$ exhibits $w$ as a product of more than $m$ and at most $k$ nonempty factors. Suppose $m$ of the factors in this product are distinguished. If not all these factors are letters $a_i$, distinguish more letters to bring the total of distinguished letters $a_i$ to $m$. By definition of $Z$ there is a word $u \in Z$ such that $A \overset{*}{\to} u \precneqq a_1 \cdots a_t$ and $u$ contains all the distinguished letters of $a_1 \cdots a_t$. It follows that $v = w'uw''$ contains the distinguished factors of $w$ and satisfies all the conditions of Theorem A.

Finally $\lambda(p) = Af\omega$ implies $\beta(p) = z_1 \cdots z_t$ with $m < t \leq C^2$ and each $z_i \in N \cup \Sigma$. Further $\gamma(p) = \beta(p) \cdot \omega$. Consequently $w = w'u_1 \cdots u_t w''$ where each $u_i$ is the subword derived from $z_i \cdot \omega$ in the derivation $S \overset{*}{\to} w$. Because $G$ is in normal form, none of the $u_i$'s is the empty word. As before there exists $\alpha \in Z$ such that $Af \overset{*}{\to} \alpha \precneqq \beta(p)$ and $\alpha$ contains all the $z_i$'s for which $u_i$ is distinguished. We have $\alpha \cdot \omega \overset{*}{\to} u$ where $u$ is the subproduct of $u_1 \cdots u_t$ corresponding to the $z_i$'s in $\alpha$. It follows that $v = w'uw''$ satisfies the conditions of Theorem A.

DEPARTMENT OF MATHEMATICS, STEVENS INSTITUTE OF TECHNOLOGY, HOBOKEN, NJ 07030

*E-mail address*: rgilman@apollo.stevens-tech.edu